\documentclass{amsart}
\usepackage{amsfonts,amssymb,amsmath,amsthm,adjustbox,tikz-cd}
\usepackage[hidelinks]{hyperref}
\usepackage{enumerate}
\usepackage{microtype}
\usepackage{mathtools}
\usepackage{cite}
\usepackage{color}
\usepackage{float}
\usepackage{multicol}
\usepackage{longtable}
\usepackage[OT2, T1]{fontenc}
\usepackage{verbatim}

\def\Z{\mathbb{Z}}
\def\Q{\mathbb{Q}}
\def\Qbar{\overline{\Q}}
\def\F{\mathbb{F}}
\def\P{\mathbb{P}}
\def\R{\mathbb{R}}
\def\C{\mathbb{C}}
\def\Gal{\operatorname{Gal}}
\def\GL{\operatorname{GL}}
\def\SL{\operatorname{SL}}

\def\det{\operatorname{det}}

\def\im{\operatorname{im}}

\def\Aut{\operatorname{Aut}}

\def\ord{\operatorname{ord}}

\newcommand{\mc}[1]{\href{https://beta.lmfdb.org/ModularCurve/Q/#1/}{\texttt{#1}}}

\theoremstyle{plain}
\newtheorem{theorem}{Theorem}

\newtheorem{proposition}[theorem]{Proposition}
\newtheorem{remark}[theorem]{Remark}
\newtheorem{step}{Step}

\theoremstyle{definition}

\title[Rational points on modular curves via maps to elliptic curves]{Rational points on modular curves via maps to elliptic curves with rank zero}
\date{\today}
\author{Jacob Mayle}
\address{Jacob Mayle, Department of Mathematical Sciences, University of Delaware, Newark, DE 19716}
\email{mayle@udel.edu}
\author{Jeremy Rouse$^{\ast}$}
\address{Jeremy Rouse, Department of Mathematics, Wake Forest University, Winston-Salem, NC 27109}
\email{rouseja@wfu.edu}
\thanks{$^\ast$Correspondence: Jeremy Rouse, Wake Forest University, rouseja@wfu.edu}

\begin{document}

\begin{abstract}
A fundamental problem in arithmetic geometry is to determine the image of the mod $N$ Galois representation for all elliptic curves over $\mathbb{Q}$ and integers $N \geq 1$. For a given subgroup $G \le \mathrm{GL}_2(\mathbb{Z}/N\mathbb{Z})$, there is a modular curve $X_G$ whose rational points parametrize elliptic curves for which the image of the mod $N$ Galois representation is contained in $G$. If $X_G$ admits a map to an elliptic curve $E/\mathbb{Q}$ for which $E(\mathbb{Q})$ has rank $0$, then its rational points can be effectively determined, provided that a map $X_G \to E$ is known. In this article, we give a method for constructing such maps. Using this method, together with existing methods and results, we systematically determine the rational points of $X_G$ for more than $99\%$ of modular curves of level at most $70$.
\end{abstract}

\maketitle

\section{Introduction}
\label{section:intro}

If $E$ is an elliptic curve defined over a number field $K$ without complex multiplication, then Serre's open image theorem \cite{MR387283} implies that the image of the adelic Galois representation
\(
\rho_{E} \colon \Gal(\overline{K}/K) \to \GL_{2}(\widehat{\Z})
\)
has finite index in $\GL_{2}(\widehat{\Z})$, and hence coincides with the full preimage of the image of the mod $N$ Galois representation
\(
\rho_{E,N} \colon \Gal(\overline{K}/K) \to \GL_{2}(\Z/N\Z)
\)
for some integer $N \geq 1$. The properties of these Galois representations have played a crucial role in many major developments in modern number theory, including the proof of Fermat's Last Theorem \cite{Wiles}, the inverse Galois problem for ${\rm PSL}_{2}(\F_{p})$ (see, for example, \cite{Shih,ZywinaModular,ZywinaInverse}), and the determination of all imaginary quadratic fields of class number one (see \cite{Baker, Heegner, Stark}).

This article is concerned with Mazur's Program B \cite[p.\ 109]{Mazur}, which asks the following: given a number field $K$ and a subgroup $G \leq \GL_{2}(\widehat{\Z})$, classify all elliptic curves $E/K$ for which the image of $\rho_{E}$ is contained in $G$. The first case to consider is $K = \Q$, which is the focus of this article. Given a subgroup $G \leq \GL_{2}(\Z/N\Z)$ with surjective determinant, the associated modular curve $X_{G}$ is a smooth, projective, geometrically integral curve
defined over $\Q$. The rational points on $X_{G}$ parametrize elliptic curves $E/\Q$ for which $\im \rho_{E,N} \leq G$ (Proposition~\ref{P:ModCrvsParam}). Motivated by Mazur's Program B,
a community of mathematicians has contributed to building a database of modular curves in the $L$-functions and Modular Forms Database (LMFDB) \cite{lmfdb}. A beta version of this database is available at \url{https://beta.lmfdb.org/ModularCurve/Q/}.
At present, the database contains an entry for every subgroup $G$ with surjective determinant and level $N \leq 70$, including data on rational cusps, rational CM points, the decomposition of the Jacobian of $X_{G}$ up to isogeny, and in many cases a model of $X_{G}$ together with maps from $X_{G}$ to other
modular curves.

Let $X$ be a smooth projective curve of genus at least $2$. By Faltings's theorem, the set of rational points $X(\Q)$ is finite, but the theorem provides no general method for computing this set. One approach is to exhibit a nonconstant morphism (which we simply call a \emph{map}) $\phi \colon X \to Y$ to a curve $Y$ whose set of rational points $Y(\Q)$ is finite and effectively computable. In this situation, 
\[
X(\Q) \subseteq \bigcup_{P \in Y(\Q)} \phi^{-1}(P),
\]
and the right-hand side is finite, allowing $X(\Q)$ to be computed effectively. The goal of this article is to develop and apply this strategy in the context of modular curves over $\Q$ that admit a map to an elliptic curve over $\Q$ of rank $0$ (in fact, our technique for constructing such a map is independent of the rank of the elliptic curve).

The LMFDB lists all $804450$ open subgroups of $\GL_{2}(\widehat{\Z})$ up to conjugacy of level $\leq 70$, with surjective determinant and containing $-I$. These correspond to modular curves $X_{G}$ defined over $\Q$ whose rational points dictate $\im \rho_{E,N}$ for $N \leq 70$. For these curves, the LMFDB includes the factorization of $J(X_{G})$ up to isogeny as a product of modular abelian varieties. By counting those for which the decomposition of \(J(X_{G})\) contains a rank \(0\) elliptic curve factor, we find that there are $797591$ subgroups $G$ such that the modular curve $X_{G}$ admits a map to an elliptic curve over $\Q$ of rank $0$. (The high proportion is due in part to covering relations: if \(X_G\)
admits a map to a rank \(0\) elliptic curve \(E/\Q\) and \(H \subseteq G\),
then the natural map \(X_H \to X_G\) gives, by composition, a map
\(X_H \to E\).) Of these $797591$ modular curves, we are able to rigorously determine the rational points of all but ten in this article. For some of these curves (for example $X_{0}(N)$, or those of $2$-power level), the rational points were already known. We list the remaining $6859 \, (= 804450 - 797591)$ modular curves for which our method does not apply in our GitHub repository \cite{ModCrvToEC}. Among these curves, $2137$ have genus $0$ and $1401$ have genus $1$ and positive rank; the remaining $3321$ curves have only finitely many rational points. Many of these curves could be handled using known methods (such as local methods, Chabauty, and \'etale descent), though some are beyond reach at the moment. 

Throughout this article, we refer to modular curves using the labels assigned in the LMFDB (see Section~\ref{section:back} for a description of these labels). Our main result is the following.

\begin{theorem}
\label{theorem:main}
Let $X_{G}$ be a modular curve over $\Q$ with level $\leq 70$, where $G$ contains $-I$ and has surjective determinant, and suppose that $X_G$ admits a map to an elliptic curve over $\Q$ of rank $0$. If $X_{G}$ is not among the following ten curves
\begin{multicols}{2}
\begin{itemize}
\item \mc{58.812.61.a.1}
\item \mc{58.812.61.b.1}
\item \mc{62.930.63.a.1}
\item \mc{57.1026.67.a.1}
\item \mc{57.1026.72.c.1}
\item \mc{46.1518.104.a.1}
\item \mc{69.1518.104.a.1}
\item \mc{69.1518.115.b.1}
\item \mc{58.2436.175.a.1}
\item \mc{55.3300.239.a.1}
\end{itemize}
\end{multicols}
\noindent then every point in $X_{G}(\Q)$ is either a cusp, a CM point, or has image under the $j$-map appearing in Table~\ref{jinvtable}.
\end{theorem}

\begin{remark}
In particular, we provably determine the rational points on every modular curve $X_{G}$ of level $\leq 70$ that maps to a rank $0$ elliptic curve over $\Q$ for which the genus of $X_{G}$ is $\leq 60$. In addition, each of the ten modular curves listed in Theorem~\ref{theorem:main}, with the exception of \mc{55.3300.239.a.1}, maps
to some modular curve $X_{{\rm ns}}^{+}(\ell)$ for $\ell \in \{ 19, 23, 29, 31 \}$. The rational points on these curves are unknown and these cases are expected to be difficult.
\end{remark}

\begin{remark}
Our method for computing a map $X_{G} \to E$ does not rely on the elliptic curve $E/\Q$ having rank $0$, and thus our code can also be used to produce
maps to elliptic curves of positive rank. One potential application is that if $X_{G}$ admits two independent maps to the same rank $1$ elliptic curve,
then $X_{G}(\Q)$ can be determined using the method of Dem'janenko--Manin. For more details on this method, see \cite[Subsection 13.3.1]{CohenVolII}.
\end{remark}

\begin{remark}
While the LMFDB does not contain the decomposition of the Jacobian of $X_{G}$ when the level of $G$ is $> 70$, our code contains an intrinsic to search for elliptic curves over $\Q$ that might occur as isogeny factors of $J(X_{G})$. Once a candidate elliptic curve factor $E/\Q$ of $J(X_{G})$ has been identified, our code can construct a map from $X_G$ to an elliptic curve isogenous to $E$ over $\Q$, provided that $E$ is indeed a factor of the Jacobian up to isogeny.
\end{remark}

For our computations, we need only consider modular curves $X_G$ that are \emph{minimal} in the following sense:
\begin{enumerate}
    \item \label{Pt1} The curve $X_G$ admits a map to an elliptic curve over $\Q$ of rank $0$, and
    \item There is no open subgroup $\widetilde{G} \leq \GL_2(\widehat{\Z})$ properly containing $G$ for which $X_{\widetilde{G}}$ satisfies \eqref{Pt1}.
\end{enumerate}
There are exactly $1034$ minimal modular curves of level $\leq 70$. Once the rational points on a minimal modular curve have been determined, we can check which of these lift to any given modular curve that covers the minimal one (this can be accomplished using either \cite{Novak} or \cite{ZywinaOpen}). The table below lists all non-CM $j$-invariants arising from rational points on any of the $797581 \, (= 797591 - 10)$ modular curves for which we provably determined the rational points. For each such $j$-invariant, the table also records the label(s) of the corresponding minimal modular curve(s), an example of an elliptic curve $E_j$ with that $j$-invariant, and the adelic level and adelic index of $E_j$. All of the $j$-invariants appearing in the table have previously been identified by Zywina \cite{ZywinaOpen, ZywinaModularGitHub}.

{\setlength{\tabcolsep}{3pt}
\scriptsize 
\begin{longtable}{|l|l|l|r|r|} \caption{Non-cuspidal, non-CM $j$-invariants arising from minimal $X_{G}$} \\
\hline $j$-invariant & Label of $X_{G}$ & Elliptic curve & Adelic level & Adelic index \\ \hline 
\endhead
\hline \multicolumn{5}{|r|}{\textit{Continued on next page}} \\ \hline
\endfoot
\hline
\endlastfoot
\hline \label{jinvtable}
$-2^{9} \cdot 3^{3}$ & \mc{6.6.1.b.1} & $[ 0, 0, 0, -24, -48 ]$ & 12 & 24\\
$2^{-7} \cdot 3^{3} \cdot 5 \cdot 7^{5}$ & \mc{7.56.1.b.1} & $[ 1, -1, 0, -107, -379 ]$ & 280 & 224\\
$-11 \cdot 131^{3}$ & \mc{11.12.1.a.1} & $[ 1, 1, 1, -30, -76 ]$ & 88 & 480\\
$-11^{2}$ & \mc{11.12.1.a.1} & $[ 1, 1, 0, -2, -7 ]$ & 88 & 480\\
$-2^{-2} \cdot 3^{3} \cdot 11^{3}$ & \mc{12.16.1.a.1} & $[ 1, -1, 0, -6, 8 ]$ & 12 & 128\\
$2^{-6} \cdot 3^{2} \cdot 23^{3}$ & \mc{12.16.1.a.1} & $[ 1, -1, 1, 4, -1 ]$ & 12 & 128\\
$-2^{10} \cdot 3^{4}$ & \mc{12.24.1.g.1} & $[ 0, 0, 0, -36, 84 ]$ & 12 & 24\\
$-2^{10} \cdot 3^{3} \cdot 5$ & \mc{12.24.1.h.1} & $[ 0, 0, 0, -60, 180 ]$ & 12 & 24\\
$-2^{6} \cdot 3^{3} \cdot 23^{3}$ & \mc{12.24.1.n.1} & $[ 0, 0, 0, -138, -624 ]$ & 72 & 144\\
$2^{-6} \cdot 3^{3} \cdot 5^{3}$ & \mc{12.24.1.o.1} & $[ 1, -1, 1, 10, 69 ]$ & 132 & 96\\
$2^{15} \cdot 7^{5}$ & \mc{14.42.1.a.1} & $[ 0, -1, 1, -45047, -3665012 ]$ & 2758 & 84\\
$-2^{4} \cdot 3^{2} \cdot 13^{3}$ & \mc{15.20.1.a.1} & $[ 0, 0, 0, -39, 94 ]$ & 120 & 160\\
$2^{4} \cdot 3^{3}$ & \mc{15.20.1.a.1} & $[ 0, 0, 0, 9, -18 ]$ & 120 & 160\\
$-2^{-3} \cdot 5^{2} \cdot 241^{3}$ & \mc{15.24.1.a.1} & $[ 1, 0, 1, -126, -552 ]$ & 120 & 384\\
$-2^{-5} \cdot 5 \cdot 29^{3}$ & \mc{15.24.1.a.1} & $[ 1, 1, 1, -3, 1 ]$ & 120 & 384\\
$-2^{-1} \cdot 5^{2}$ & \mc{15.24.1.a.1} & $[ 1, 0, 1, -1, -2 ]$ & 120 & 384\\
$2^{-15} \cdot 5 \cdot 211^{3}$ & \mc{15.24.1.a.1} & $[ 1, 1, 1, 22, -9 ]$ & 120 & 384\\
$-2^{-15} \cdot 29^{3} \cdot 41^{3}$ & \mc{15.36.1.b.1} & $[ 1, 1, 1, -322, 2127 ]$ & 1560 & 576\\
$2^{-3} \cdot 11^{3}$ & \mc{15.36.1.b.1} & $[ 1, 1, 1, 3, -5 ]$ & 1560 & 576\\
$-2^{-1} \cdot 17^{2} \cdot 101^{3}$ & \mc{17.18.1.a.1} & $[ 1, 0, 1, -3041, 64278 ]$ & 680 & 576\\
$-2^{-17} \cdot 17 \cdot 373^{3}$ & \mc{17.18.1.a.1} & $[ 1, 1, 0, -660, -7600 ]$ & 680 & 576\\
$2^{12} \cdot 3^{3}$ & \mc{18.18.1.b.1} & $[ 0, 0, 1, -21, -37 ]$ & 126 & 36\\
$-2^{12} \cdot 3^{7} \cdot 5$ & \mc{18.54.2.b.1} & $[ 0, 0, 1, -135, -604 ]$ & 18 & 54\\
$-2^{-2} \cdot 5^{2} \cdot 41^{3}$ & \mc{20.24.1.g.1} & $[ 1, 1, 1, -363, -2819 ]$ & 340 & 192\\
$2^{-10} \cdot 5 \cdot 59^{3}$ & \mc{20.24.1.g.1} & $[ 1, 0, 1, 104, 358 ]$ & 340 & 192\\
$-2^{-10} \cdot 3^{3} \cdot 37^{3} \cdot 47$ & \mc{20.40.2.b.1} & $[ 1, -1, 1, -109, 469 ]$ & 940 & 80\\
$-2^{2} \cdot 3^{2}$ & \mc{20.40.2.b.1} & $[ 0, 0, 0, -3, 14 ]$ & 20 & 80\\
$-2^{-10} \cdot 3^{3} \cdot 5^{4} \cdot 11^{3} \cdot 17^{3}$ & \mc{20.60.3.w.1} & $[ 1, -1, 0, -51717, 4539861 ]$ & 3540 & 240\\
$2^{16} \cdot 3^{3} \cdot 5^{-5} \cdot 17$ & \mc{20.120.6.a.1} & $[ 0, 0, 0, -272, -1564 ]$ & 20 & 120\\
$2^{10} \cdot 3^{2} \cdot 5^{-5} \cdot 79^{3}$ & \mc{20.120.6.a.1} & $[ 0, 0, 0, -948, -11228 ]$ & 20 & 120\\
$-2^{-7} \cdot 3^{3} \cdot 5^{3} \cdot 383^{3}$ & \mc{21.32.1.a.1} & $[ 1, -1, 0, -1077, 13877 ]$ & 504 & 768\\
$-2^{-3} \cdot 3^{2} \cdot 5^{6}$ & \mc{21.32.1.a.1} & $[ 1, -1, 1, -5, 5 ]$ & 504 & 768\\
$-2^{-21} \cdot 3^{2} \cdot 5^{3} \cdot 101^{3}$ & \mc{21.32.1.a.1} & $[ 1, -1, 1, -95, -697 ]$ & 504 & 768\\
$2^{-1} \cdot 3^{3} \cdot 5^{3}$ & \mc{21.32.1.a.1} & $[ 1, -1, 0, 3, -1 ]$ & 504 & 768\\
$-2^{6} \cdot 3^{3} \cdot 23^{3}$ & \mc{24.24.2.d.1} & $[ 0, 0, 0, -138, -624 ]$ & 72 & 144\\
$2^{6} \cdot 3^{-6} \cdot 31^{3}$ & \mc{24.72.2.a.1} & $[ 0, -1, 0, -289, -1007 ]$ & 168 & 144\\
$2^{6} \cdot 3^{-3} \cdot 109^{3}$ & \mc{24.72.2.f.1} & $[ 0, 1, 0, -254, 1476 ]$ & 168 & 144\\
$-2^{-14} \cdot 3^{3} \cdot 13 \cdot 479^{3}$ & \mc{28.32.2.a.1} & $[ 1, -1, 0, -389, -2859 ]$ & 364 & 768\\
$2^{-2} \cdot 3^{3} \cdot 13$ & \mc{28.32.2.a.1} & $[ 1, -1, 0, 1, 1 ]$ & 364 & 768\\
$-2^{12} \cdot 3^{-1} \cdot 5^{2}$ & \mc{30.12.1.d.1} & $[ 0, -1, 1, -8, -7 ]$ & 30 & 48\\
$2^{12} \cdot 3^{-5} \cdot 5$ & \mc{30.12.1.d.1} & $[ 0, 1, 1, 2, 4 ]$ & 30 & 48\\
$-2^{9} \cdot 3^{-5} \cdot 5^{4} \cdot 11^{3}$ & \mc{30.30.2.b.1} & $[ 0, -1, 0, -733, -7403 ]$ & 60 & 60\\
$2^{12}$ & \mc{30.36.1.q.1} & $[ 0, 1, 1, -12, -17 ]$ & 1110 & 288\\
$2^{12} \cdot 211^{3}$ & \mc{30.36.1.q.1} & $[ 0, 1, 1, -2602, 50229 ]$ & 1110 & 288\\
$-3^{3} \cdot 5^{4} \cdot 11 \cdot 19^{3}$ & \mc{30.60.4.a.1} & $[ 1, -1, 1, -980, -11558 ]$ & 30 & 120\\
$-2^{2} \cdot 3^{7} \cdot 5^{3} \cdot 439^{3}$ & \mc{36.108.6.g.1} & $[ 0, 0, 0, -1126035, 459913278 ]$ & 684 & 216\\
$-7 \cdot 137^{3} \cdot 2083^{3}$ & \mc{37.38.2.a.1} & $[ 1, 1, 1, -208083, -36621194 ]$ & 5180 & 2736\\
$-7 \cdot 11^{3}$ & \mc{37.38.2.a.1} & $[ 1, 1, 1, -8, 6 ]$ & 5180 & 2736\\
$-2^{3} \cdot 5^{4}$ & \mc{40.30.2.e.1} & $[ 0, 1, 0, -8, 8 ]$ & 40 & 60\\
$-2^{6} \cdot 719^{3}$ & \mc{45.54.2.c.1} & $[ 0, 1, 0, -6231, 187253 ]$ & 9360 & 864\\
$2^{6}$ & \mc{45.54.2.c.1} & $[ 0, 1, 0, 9, 53 ]$ & 9360 & 864\\
$17^{3}$ & \mc{48.72.3.x.1} & $[ 1, 1, 1, -23, -44 ]$ & 3120 & 1152\\
$257^{3}$ & \mc{48.72.3.x.1} & $[ 1, 1, 1, -348, -2644 ]$ & 3120 & 1152\\ \hline
\end{longtable} }

The problem of finding maps from the modular curves $X_{0}(N)$ to elliptic curves was considered in detail by Cremona in \cite[Chapter 2]{Cremona-alg}. In comparison with the setting of $X_0(N)$, the general situation for an arbitrary modular curve $X_G$ is complicated by several factors. These include the fact that $X_G$ need not have a rational cusp that can serve as a base point for the Abel--Jacobi map, the Fourier expansions of modular forms for $X_G$ need not be expressible in terms of integer powers of $q = e^{2 \pi i z}$, and the Fourier coefficients need not be rational. We make extensive use of Zywina's algorithms for computing modular forms for arbitrary subgroups $G \leq \GL_{2}(\Z/N\Z)$, as described in \cite{ZywinaModular} and implemented in his repository \cite{ZywinaModularGitHub}. Further, since our goal is to compute the rational points on $X_G$, we must also explicitly construct the map $X_G \to E$ in terms of a defining equation for $X_G$ and use this map to determine all rational points on $X_G$.

Following a review of relevant background in Section~\ref{section:back}, we describe our algorithm in detail in Section~\ref{section:alg}, and in Section~\ref{section:ex} we illustrate the method by working through a numerical example. In Section~\ref{section:results}, we discuss the details of our computations. In particular, we successfully ran our code on all minimal modular curves of genus $\leq 10$ that have at least one rational point. The runtime of our algorithm generally increases with the genus, and the highest genus curve for which we successfully ran our code was \mc{63.756.54.a.1} of genus $54$, which required approximately $9$ hours of computation.

In Section~\ref{section:unique_cases}, we discuss several unique cases that arise, the most notable of which is the high genus modular curve \mc{63.1512.115.t.1}, which covers the genus $3$ curve $X$ with LMFDB label \mc{21.84.3.a.1}. Although this curve $X$ does not admit a map to a rank $0$ elliptic curve over $\Q$, it does admit a map to a rank $0$ elliptic curve over $\Q(\sqrt{-3})$. We adapt our approach to compute this map and determine the set $X(\Q(\sqrt{-3}))$. This technique should be applicable to modular curves
$X_{G}$ whose Jacobians have a factor $A$ that is the restriction of scalars to $\Q$ of
a $\Q$-curve with rank $0$, but we do not attempt a general implementation. There is also the possibility of using this technique combined with elliptic curve Chabauty in situations where $J(X_{G})$ has a simple factor $A$ with $\dim(A) = {\rm rank}~A(\Q)$ which occurs as a factor of ${\rm Res}_{K/\Q}(E)$ for some elliptic curve $E$ for which $[K : \Q] > \dim(A)$. Finally, in Section~\ref{section:remaining_cases} we complete the proof of Theorem~\ref{theorem:main} by handling minimal covers of the ten modular curves listed in the theorem.

\subsection{Acknowledgments}

We thank David Zywina and the anonymous reviewers for many helpful comments.

\section{Background}
\label{section:back}

Let $K$ be a number field and $E/K$ be an elliptic curve. Write $E[N]$ for the subgroup of $E(\overline{K})$ consisting of points of order dividing $N$. The natural action of $\Gal(\overline{K}/K)$ on $E(\overline{K})$ restricts to an action on $E[N]$, giving rise to a homomorphism
\( \Gal(\overline{K}/K) \to \Aut(E[N]) \)
defined by 
\(
\sigma \mapsto \sigma|_{E[N]}.
\)
Since $E[N] \simeq (\Z/N\Z)^{2}$, a choice of basis allows us
to identify $\Aut(E[N])$ with $\GL_{2}(\Z/N\Z)$. We regard elements of $E[N]$ as column vectors, with $\GL_2(\Z/N\Z)$ acting by left multiplication; this convention agrees with that used in the LMFDB. The resulting homomorphism
\[ \rho_{E, N} \colon \Gal(\overline{K}/K) \longrightarrow \GL_{2}(\Z/N\Z) \]
is called the \emph{mod $N$ Galois representation of $E$}. Since this map depends on the choice of basis, its image can only be considered up to conjugation in $\GL_{2}(\Z/N\Z)$. Taking the inverse limit over all integers $N \geq 1$ (with a compatible choice of bases) yields the \emph{adelic Galois representation of $E$},
\[
\rho_E \colon \Gal(\overline{K}/K) \longrightarrow \GL_2(\widehat{\Z})
\]
which describes the action of $\Gal(\overline{K}/K)$ on the adelic Tate module of $E$. Here $\widehat{\Z} = \varprojlim \Z/N\Z$ denotes the ring of profinite integers. As before, the image of $\rho_E$ can only be considered up to conjugation in $\GL_{2}(\widehat{\Z})$.

Let $G \leq \GL_2(\widehat{\Z})$ be an open subgroup. Let $N$ be the \emph{level} of $G$, that is, the smallest positive integer such that $G$ contains the kernel of the reduction map $\GL_2(\widehat{\Z}) \to \GL_2(\Z/N\Z)$. We often identify $G$ with its reduction in $\GL_2(\Z/N\Z)$. Associated to $G$ is a modular curve $X_{G}$ whose points parametrize elliptic curves with a $G$-level structure (for more details, see \cite[Section IV-3]{DeligneRapoport} and \cite{Zywina2015}). The set $X_{G}(\C)$
can be identified with the compactification of the upper half-plane $\mathbb{H}$ modulo $\Gamma = G \cap \SL_{2}(\Z)$. We assume throughout that $\det \colon G \to (\Z/N\Z)^{\times}$ is surjective. When this holds, $X_{G}$ is defined over $\Q$ and is smooth, projective, and geometrically integral. The space of holomorphic $1$-forms on $X_{G}$ can be identified
with a $\Q$-vector space of weight $2$ cusp forms for $\Gamma$. If $G \leq \widetilde{G}$, there is a natural morphism $X_{G} \to X_{\widetilde{G}}$. An important instance of this is the \emph{$j$-map} $j \colon X_{G} \to \P^{1}$, which is obtained when $\widetilde{G} = \GL_{2}(\Z/N\Z)$. A \emph{cusp} of $X_{G}$ is a point on $X_{G}$ whose image under $j$ is $(1 : 0)$. The term \emph{modular curve} is sometimes used more broadly than we do here, such as in \cite{MR879929} to include Atkin--Lehner quotients as well, but throughout this article we only consider modular curves of the form \(X_G\) for subgroups \(G \leq \GL_2(\widehat{\mathbb Z})\) satisfying the hypotheses above.

In the spirit of Mazur's Program B, we are motivated to study modular curves because the $K$-rational points on $X_{G}$ parametrize elliptic curves $E$ whose mod $N$ Galois representation has image contained in $G$, as made precise by the proposition below.

\begin{proposition} \cite[Proposition 3.3]{Zywina2015} \label{P:ModCrvsParam}
Let $G \leq \GL_{2}(\Z/N\Z)$ be a subgroup that contains $-I$ and has surjective determinant. Suppose that $E/\Q$ is an elliptic curve with $j(E) \not\in \{0, 1728\}$. Then $\im \rho_{E,N}$ is conjugate to a subgroup of $G$ if and only if
$j(E) \in j(X_{G}(\Q))$.
\end{proposition}

\noindent When $-I \not\in G$, the modular curves
$X_{G}$ and $X_{\langle G, -I \rangle}$ coincide, and the condition $\im \rho_{E,N} \leq G$ is no longer determined solely by the $j$-invariant. In this case, there exists a universal elliptic curve $\mathcal{E}$ defined over the complement of the $j = 0$ and $j = 1728$ loci of $X_{G}$ that parametrizes 
elliptic curves whose mod $N$ Galois image is contained in $G$. For more details,
see \cite[Section 13]{ZywinaOpen}.

To compute the modular curve $X_{G}$, we rely on Zywina's code from \cite{ZywinaModular}, which follows the convention that $X_{G}$ parametrizes elliptic curves for which $\im \rho_{E,N}$ is contained in $G^{\intercal} = \{ g^{\intercal} : g \in G \}$. For this reason, we take
the transpose of our group prior to using Zywina's code to compute a model; see \cite[Remark 3.6]{ZywinaOpen}. 

An elliptic curve $E/K$ has \emph{complex multiplication} (CM) if $\operatorname{End}_{\overline{K}}(E)$ is an order in an imaginary quadratic field. If such a curve is defined over $\Q$, then $\operatorname{End}_{\Qbar}(E)$ must be an imaginary quadratic order of class number $1$, which has one of the following discriminants:
\[
-3, \, -4, \, -7, \, -8, \, -11, \, -12, \, -16, \, -19, \, -27, \, -28, \, -43, \, -67, \text{ or } -163.
\]
\noindent It is helpful to be able to determine, via group theory, the number of rational CM points on a modular curve $X_{G}$. For this purpose, we rely on the recent algorithm and results of \cite{Novak}.

When referring to modular curves, we often use the LMFDB label \texttt{N.i.g.c.n}, which records the level, index, genus, and Gassmann class, together with an integer tiebreaker distinguishing nonconjugate subgroups sharing the same first four invariants. Alternatively, we sometimes use the classical notation such as $X_0(N)$, $X_{\rm sp}(N)$, $X_{\rm ns}(N)$, $X_{\rm sp}^+(N)$, and $X_{\rm ns}^+(N)$ to denote the modular curves associated with the Borel subgroup, the split Cartan subgroup, the non-split Cartan subgroup, and the normalizers of the split and non-split Cartan subgroups, respectively. In addition, if $H_{1}$ and $H_{2}$ are two subgroups of $\GL_{2}(\widehat{\Z})$ with coprime levels that both contain $-I$, the modular curve $X_{H_{1} \cap H_{2}}$ is the smooth curve birational to the fiber product $X_{H_{1}} \times_{X(1)} X_{H_{2}}$ of $X_{H_{1}}$ and $X_{H_{2}}$ over the $j$-line.

To reduce the number of cases we need to consider, we make use of several prior
results about rational points on modular curves, especially in situations where a modular curve under consideration covers another modular curve whose rational points are already known. In particular, we use the main result of \cite{RSZB}, which applies to $\ell$-adic images for $\ell \leq 11$, and the results of
Zywina \cite[Proposition 1.13]{Zywina2015} and Furio and Lombardo \cite{FurioLombardo}, which show
that if $\ell > 5$, there is no non-CM elliptic curve $E/\Q$ for which $\im \rho_{E,\ell}$ is a proper
subgroup of the normalizer of the non-split Cartan. In addition, we use the results of
\cite{Cursed}, which show that all the rational points on $X_{{\rm sp}}^{+}(13)$ and $X_{{\rm ns}}^{+}(13)$ are cusps or CM points, and \cite[Section 5.1]{MoreQuadChab}, which determines the rational
points on the modular curve $X_{S_{4}}(13)$.

We also rely on some generalizations of the formal immersion technique of Mazur \cite{Mazuriso}. 
In \cite{Lemos1}, Lemos shows that if $E/\Q$ is a non-CM elliptic curve with
a rational cyclic isogeny and $\im \rho_{E,\ell}$ is contained in the normalizer of the non-split Cartan mod $\ell$ for some prime $\ell > 3$, then $j(E) \in \Z$. The paper then enumerates the finitely many integral $j$-invariants of elliptic curves that have a rational cyclic isogeny. Examining these cases shows that it is not possible for an elliptic curve $E/\Q$ to have a cyclic isogeny and image contained in the normalizer of the non-split Cartan mod $\ell$ for $\ell > 3$. 
The results of \cite{Lemos2} rule out the possibility of a non-CM elliptic curve $E/\Q$
having both split and non-split Cartan images for two primes $\ell$ and $p$. The result below can be proven using this method, but is strictly stronger than the results stated in \cite{Lemos2}.

\begin{theorem} \cite[Theorem 6.2]{DanielsRouse}
Let $E/\Q$ be a non-CM elliptic curve. If there is a prime $q$ for which $\im \rho_{E,q}$ is contained in the normalizer of the split Cartan mod $q$, then there is no prime $p > 3$ for which $\im \rho_{E,p}$ is contained in the normalizer of the non-split Cartan mod $p$.
\end{theorem}

We present an algorithm for computing a map $X_G \to E$ under the assumption that $X_{G}(\Q)$ is non-empty. This algorithm has several numerical and heuristic steps, and can fail in situations where these heuristics do not hold. In the latter case, rerunning the code (possibly increasing some parameters) should lead to success. Because of these heuristics, the polynomials produced are not guaranteed to define a map from $X_{G}$ to an elliptic curve $E$. To ensure the output is correct, we must verify that each map produced by the algorithm is indeed a map $X_{G} \to E$. To this end, we make use of the following well-known result about modular forms. For a proof, see \cite[Proposition 4.12]{Milne}.
\begin{proposition}
\label{checkingprop}
Suppose that $f$ is a nonzero modular form of weight $k$ for a finite index subgroup
$\Gamma$ of $\SL_{2}(\Z)$ that contains $-I$. Then
\[
  \sum_{Q \in X_{\Gamma}(\mathbb{C})} \frac{\ord_{Q}(f)}{e_{Q}}
  = k(g-1) + \frac{k \cdot \epsilon_{2}}{4}
  + \frac{k \cdot \epsilon_{3}}{3} + \frac{k \cdot \epsilon_{\infty}}{2},
\]
where $g$ is the genus of $X_{\Gamma}$ and $e_{Q}$ denotes the order of the stabilizer of $Q$ in  $\Gamma/\langle \pm I \rangle$. The quantity $\epsilon_{\infty}$ is the number of cusps of $X_{\Gamma}$, and $\epsilon_{2}$ and $\epsilon_{3}$ are the numbers of elliptic points of orders $2$ and $3$, respectively.
\end{proposition}
The Riemann--Hurwitz theorem (see \cite[Theorem 3.1.1]{DiamondShurman}) shows that the right hand side
equals $\frac{k}{12} [\SL_{2}(\Z) : \Gamma]$. As a consequence,
if $f$ is a weight $k$ modular form for $\Gamma$ and
the sum of the orders of vanishing of $f$ at all the cusps of $X_{\Gamma}$ is greater than $\frac{k}{12} [\SL_{2}(\Z) : \Gamma]$, then $f = 0$.

\section{Description of algorithm}
\label{section:alg}

In this section, we document the algorithms that are used in the paper.
The first function, {\tt InitializeModEC}, creates a Magma record that holds the various pieces of data needed for later use. This function
takes as input a level $N$, as well as matrices that generate
a subgroup $G \leq \GL_{2}(\Z/N\Z)$. It then uses Zywina's code to compute a model of the modular curve $X_{G}$ as well as Fourier expansions of modular forms invariant under the action of $G$. 

The main function is {\tt FindMapsToEC}. The input consists of a model of \(X_G\), Fourier expansions of a basis of \(S_2(G)\), and a finite list of elliptic curves over \(\Q\) (specified by their Weierstrass equations) together with their expected multiplicities as factors of the Jacobian \(J(X_G)\). The output is a map 
\( X_G \to E', \)
where \(E'/\Q\) is an elliptic curve isogenous over \(\Q\) to one of the input elliptic curves. This map is represented by homogeneous polynomials. The model of \(X_G\), the modular forms, the map $X_{G} \to E'$, and the auxiliary data are stored in a Magma record created by {\tt InitializeModEC}. It is known from \cite[Appendix A]{RSZB} that every simple factor of $J(X_{G})$ is an abelian
variety of $\GL_{2}$-type and is hence isogenous to an abelian variety $A_{f}$ for some newform $f$ via the Eichler--Shimura
correspondence \cite[Definition 6.6.3]{DiamondShurman}. The beta version of the modular curves database in the LMFDB 
provides the decomposition of $J(X_{G})$ for every subgroup $G \leq \GL_2(\widehat{\Z})$ with surjective determinant and level $N \leq 70$.

To construct a map $X_{G} \to E$, the process is to identify a weight $2$ cusp form $f$ whose Hecke eigenvalues for primes $p \nmid N$ equal $a_{p}(E)$ and whose lattice of periods $\Lambda$ matches that of $E$. We fix a base point $Q \in X_{G}(\Q)$ and let 
$\alpha$ be the natural identification $\alpha \colon X_{G}(\C) \to \mathbb{H}/(G \cap \SL_{2}(\Z))$. (Here we have assumed that a rational base point on \(X_G\) has been found through a point search; see Step~\ref{Step2} and the discussion of the optional parameter \texttt{IgnoreBase} for how we proceed when no rational point is found.)
Then, we have a map $X_{G} \to \C/\Lambda$ given by
\[
P \mapsto \int_{\alpha(Q)}^{\alpha(P)} f \, ds 
\]
and we compose this with the isomorphism $\C/\Lambda \to E(\C)$. This produces the map $X_{G} \to E$.

There are three additional auxiliary functions provided. The first such function, named \texttt{EllipticCurveQuoCandidates}, identifies potential elliptic curve
factors of the Jacobian of $X_{G}$ by using Hecke operators and returns
a list of all elliptic curves $E$, given by their Weierstrass coefficients,  that might be factors of $J(X_{G})$,
together with upper bounds on their multiplicities. The second function is {\tt RatPtsFromMaps}, which determines $X_{G}(\Q)$ using the map found by {\tt FindMapsToEC}. Finally, the function {\tt ComputeJ} produces the $j$-map $j \colon X_{G} \to \P^{1}$.

The algorithm makes various assumptions (including the existence of a rational point on $X_{G}$) and numerical computations along the way and for this reason, the map produced at the end is not guaranteed to be correct. Hence it is necessary to verify the correctness near the end (in Step 12).

We now give a detailed description of the steps in the main function.

\begin{step} \label{Step1} Compute modular forms and a model of $X_{G}$.
\end{step}

\noindent We use Zywina's code to produce a model of $X_{G}$, estimate the power series precision needed for period computations,
and ensure that the modular forms are computed to that precision. This step is handled by the function {\tt InitializeModEC}. In the case that $X_{G}$ has genus $\geq 3$ and is not geometrically hyperelliptic, this code computes the $\Q$-vector space of weight $2$ cusp forms on $X_{G}$ and uses this to compute the canonical model. In the case that $X_{G}$ is geometrically hyperelliptic, a $\Q$-rational line bundle $\mathcal{L}$ is chosen on $X_{G}$ so that $\deg(\mathcal{L}) \geq 2g+1$. We let $R$ be the graded ring 
\[ R = \bigoplus_{n \geq 0} H^{0}(X_{G}, \mathcal{L}^{\otimes n}).\]
A theorem of Mumford (see \cite[p. 55]{Mumford}) implies that $R$ is generated in degree $1$. In the non-hyperelliptic case, we let $R$
be the canonical ring and $\mathcal{L}$ the canonical bundle. The canonical ring is generated in degree $1$ by a theorem of Max Noether (see \cite[Theorem 1.6]{AS}). 

\begin{step} \label{Step2} Test local solvability and find a rational base point.    
\end{step} 

\noindent We attempt to find a rational base point $Q \in X_{G}(\Q)$. The LMFDB contains information about whether $X_{G}(\F_{q}) = \emptyset$ for primes $q \nmid N$. When the genus of $X_{G}$ is less than $30$, we search for a rational point on $X_{G}$ and test whether $X_{G}$ has points modulo $p^{k}$ for primes $p \mid N$ and for $k \leq 3$ (for $p = 3$ we test $k \leq 4$ and for $p = 2$ we test $k \leq 8$). If a prime $p$ is found for which $X_{G}(\Q_{p}) = \emptyset$, the code terminates. We perform these local checks using our function \texttt{SmartpAdic}. In the few situations where a rational point on $X_{G}$ is not found, the code attempts to proceed using
the cusp at infinity (which need not be rational) as a base point although there is no guarantee that the resulting map will be defined over $\Q$. We considered the possibility of using a rational divisor class of positive degree (which always exists) but this would make the end of Step 8 more complicated.

\begin{step} \label{Step3} Identify the correct space of weight $2$ cusp forms.
\end{step}

\noindent Steps~\ref{Step3}--\ref{Step12} are handled by the main function {\tt FindMapsToEC}.
We are given the multiplicity of each elliptic curve $E$ as a factor of $J(X_{G})$. Let $\omega$ be the invariant differential on $E$.
If $\psi \colon X_{G} \to E$ is any map, then
$\psi^{*}(\omega)$ should be a holomorphic $1$-form on $X_{G}$ that can be identified with a weight $2$ cusp form for $G$ that is an eigenform of the Hecke operator $T_{p}$ for all primes $p \nmid N$ with eigenvalue $a_{p}(E)$. (Shimura proves this in many cases in \cite[Prop. 7.19]{Shimura}.) We use the Hecke operator 
\( T(p) \begin{bsmallmatrix} 1 & 0 \\ 0 & p^{-1} \end{bsmallmatrix}^{*}\).
This operator is very similar to the operator $T(p) \begin{bsmallmatrix} p & 0 \\ 0 & 1 \end{bsmallmatrix}^{*}$ from Kato \cite{Kato}. Kato's operator commutes with the action of $\GL_{2}(\Z/N\Z)$, and our operator differs from Kato's by the scalar matrix
$\begin{bsmallmatrix} p^{-1} & 0 \\ 0 & p^{-1} \end{bsmallmatrix}$. Our Hecke operator sends a modular form
\[
\sum_{n=1}^{\infty} a(n) q^{n/w}
\]
to
\[
\sum_{n=1}^{\infty} \sigma(a(pn)) q^{n/w} + \sum_{n=1}^{\infty} b(n) q^{pn/w}
\]
if $p \nmid N$. Here $q = e^{2 \pi i z}$, $w$ is the width of the cusp, and $\sigma$ is the automorphism of $\Q(\zeta_{N})$ that sends $\zeta_{N}$ to $\zeta_{N}^{s}$ where $s$ is the multiplicative inverse of $p$ modulo $N$. While we do not compute a formula for $b(n)$, we are able in practice to determine the action of the Hecke operator on the space of weight $2$ cusp forms by considering coefficients of $q^{n/w}$ with
$p \nmid w$. We intersect kernels of \[ T(p) \begin{bmatrix} 1 & 0 \\ 0 & p^{-1} \end{bmatrix}^{*} - a_{p}(E) \]
until we arrive at a subspace with dimension equal to the multiplicity. We then find a single simple modular form $f$ in this subspace whose Fourier coefficients
at all cusps are algebraic integers. We repeat this process for each isogeny class that was input.

\begin{step} \label{Step4} Compute periods for each $f$.
\end{step}

\noindent Let $\Gamma = G \cap \SL_{2}(\Z)$. For each form $f$ in the previous step, we heuristically compute the image of the period homomorphism $\phi \colon \Gamma \to \C$ defined by
\[
\phi(g) = 2 \pi i \int_{z}^{g(z)} f(s) \, ds.
\]
The integral above does not depend on the choice of the complex number $z$, and for ease of computation we numerically evaluate the above integral for $z = i$. We compute this quantity for matrices in $\Gamma$ that are lifts to $\SL_{2}(\Z)$
of generators of $G \cap \SL_{2}(\Z/N\Z)$, as well as randomly generated matrices 
in the kernel of the reduction map $\SL_{2}(\Z) \to \SL_{2}(\Z/N\Z)$, until we have at least \(20\) matrices. The random matrices are required to have word length at most \(50\) in \(S\) and \(T\). These defaults were chosen empirically and were typically sufficient in our computations to recover the period lattice (they can be modified by changing the parameters \texttt{NumMats} and \texttt{MaxLen} in
\texttt{RandWords}). We can compute the integral above by writing the image of $g$ in $\operatorname{PSL}_2(\Z/N\Z)$ as a word in $S = \begin{bsmallmatrix} 0 & -1 \\ 1 & 0 \end{bsmallmatrix}$ and $T = \begin{bsmallmatrix} 1 & 1 \\ 0 & 1 \end{bsmallmatrix}$, namely
\[
g = \pm  w_{1} w_{2} \cdots w_{k},
\]
where each $w_{i}$ is a power of either $S$ or $T$. Then
\[
2 \pi i \int_{i}^{g(i)} f(s) \, ds
= \sum_{j=1}^{k} 2 \pi i \int_{w_{1} w_{2} \cdots w_{j-1}(i)}^{w_{1} w_{2} \cdots w_{j}(i)} f(s) \, ds.
\]
We then use the transformation law
\begin{equation} \label{E:transformation}
\int_{w(\alpha)}^{w(\beta)} f(s) \, ds = \int_{\alpha}^{\beta} (f|w)(s) \, ds
\end{equation}
with $w = w_{1} w_{2} \cdots w_{j-1}$. The new upper limit of integration is $i-1$, $i$, or $i+1$ and so we can compute an antiderivative
of $f|w$ and evaluate this numerically at a complex number with imaginary part equal to $1$, which ensures fast convergence. The availability of Fourier expansions at all cusps greatly simplifies these period computations, in contrast to the method used by Cremona for
$\Gamma_{0}(N)$ in \cite[Sections 2.10-2.12]{Cremona-alg}. The power series precision chosen at the start of the computation ensures that these period computations are accurate to the desired amount of decimal precision.

\begin{step} \label{Step5} Compute the period lattice for each $f$. \end{step}

\noindent The previous step computes 20 periods of the form $\phi(g)$. We use that the period lattice for an elliptic curve over $\Q$ has the shape $\{ a \omega_{1} + b \omega_{2} : a, b \in \Z \}$ where $\omega_{1} \in \R$ and either $\omega_{2}$ is strictly imaginary, or ${\rm Re}(\omega_{2}) = \frac{1}{2} \omega_{1}$. We then identify the minimal nonzero real part of a period and the minimal nonzero imaginary part of a period. We divide all the real parts by this minimal nonzero real period and assume that the result is a rational number with numerator and denominator $\leq 100$. (Because of how the matrices are chosen we expect that Step~\ref{Step4} will generate simple elements of the period lattice, so these numerator and denominator bounds should be sufficient.) We do the same with the imaginary parts. Under these assumptions, we are able to determine the lattice spanned by the periods we have computed, which we assume is the image of $\phi$.

Note that this step is heuristic: it is possible for the finite set of periods computed in Step~\ref{Step4} to determine only a proper sublattice of the full period lattice. If such a failure occurs, it is likely to be detected in Step~\ref{Step6}, where we compare with the period lattices of the elliptic curves in the relevant isogeny class. If the computed lattice is a proper homothetic sublattice of the correct one, the algorithm may still produce a valid map, but one composed with multiplication by an integer on the elliptic curve and hence not of minimal degree. In any case, the resulting map is certified in Step~\ref{Step12}.

\begin{step} \label{Step6} For each modular form $f$, determine the optimal elliptic curve in each isogeny class.
\end{step} 

\noindent In the next step, we determine, for each modular form $f$, a unique elliptic curve $E'$ in the isogeny class and the Manin constant $c$. This is the proportionality constant between $f$ and the pullback of the invariant differential on $E'$. (In a situation where $J(X_{G})$ admits multiple independent maps to $E'$, the dimension of the space of modular forms with Hecke eigenvalues matching those of $E'$ is equal to the number of independent maps by the results of Appendix A of \cite{RSZB} and so we expect that every such $f$ is the pullback of the invariant differential of $E'$ under some map $X_{G} \to E'$.) To determine the optimal curve in the isogeny class, we divide the two periods for $f$ by those of each elliptic curve in the isogeny class until we find one where the ratios are approximately the same. We then determine the Manin constant $c$. In the case of $X_{0}(N)$, it is known \cite{Manin1} that $c$ is an integer, and for an optimal curve it was conjectured by Manin that $c = 1$. This conjecture is known (see \cite{Manin2}) when $E$ is semistable, and
in \cite{Manin3} it is shown that this constant divides the degree of the map $X_{0}(N) \to E$ (apart from factors of $2$ and $3$ that could appear in certain cases).

In our case, $c$ need not be an integer, but in all of the examples we have
computed we have observed that the numerator of $c$ is $1$ and that the
primes dividing the denominator of $c$ also divide $N$. To compute $c$ exactly from a numerical approximation, we make the assumption that
the numerator of $c$ is $1$. The result is a modular form $cf$ with the property that
the integrals $\int_{\infty}^{P} cf \, ds$ are well-defined modulo $\Lambda(E')$, the period lattice for $E'$. 

\begin{step} \label{Step7} Compute degrees of the maps to the elliptic curves and pick the one with smallest degree.
\end{step} 

\noindent Next, for each candidate modular form $f$, we use \cite[Theorem 3]{Cremona} to compute the degree $r$ of the corresponding map $X_{G} \to E'$. We then pick the pair $(f,E')$ for which the corresponding $r$ is minimal. In some cases, this is the slowest part of the algorithm, because computing each degree requires the computation of $[\SL_{2}(\Z) : \Gamma]$ periods, and there can be many options for $E'$. The Riemann--Roch theorem
ensures that we can obtain the $x$ and $y$ modular functions as ratios of elements of degree $d$ in the graded ring $R$, where $d$ satisfies
\[
\left\lfloor \frac{3r+{\rm genus}(X_{G})-1}{\deg \mathcal{L}} \right\rfloor \leq d \leq \frac{3r+{\rm genus}(X_{G})-1}{\deg \mathcal{L}} + 1.
\]
At this stage, we also ensure that enough coefficients of modular forms for $X_{G}$ have been computed to use Proposition~\ref{checkingprop} to verify our map $X_{G} \to E'$ in Step~\ref{Step12}.

\begin{step} \label{Step8} Compute the Fourier expansion of $\int_{\infty}^{P} cf \, ds$ at all cusps. \end{step}

\noindent This only requires knowing the Fourier expansion of $cf$ and the values of the integrals $\int_{\infty}^{g(\infty)} cf \, ds$.
A theorem of Manin \cite{Manin} and Drinfeld \cite{Drinfeld} implies that the difference of any two cusps on a modular curve is a torsion point,
and so the integrals above are rational linear combinations of the basis periods for $E'$. We compute these integrals numerically using techniques similar to those employed in Step~\ref{Step4}.
We also identify the images of the corresponding points in $E'(\Q(\zeta_{N}))$ by factoring division polynomials for $E'$ over $\Q(\zeta_{N})$.

\begin{step} \label{Step9} Compute Fourier expansions of the $x$ and $y$ modular functions. \end{step}

\noindent The next task is to compute Fourier expansions of the modular functions
$x : X_{G} \to E' \to \P^{1}$ and $y : X_{G} \to E' \to \P^{1}$ at all cusps. The process for computing these
is to compose the maps $\phi_{1}$ and $\phi_{2}$. Here $\phi_{1} \colon X_{G} \to \C/\Lambda(E')$ is given by
$\phi_{1}(P) = \int_{\infty}^{P} cf \, ds$, where $cf$ is the weight $2$ cusp form whose period lattice is $\Lambda(E')$.
The map $\phi_{2}$ is the isomorphism $\C/\Lambda(E') \simeq E'(\C)$. This requires knowing the Laurent series expansions at $z = 0$
of $\wp(z)$ and $\wp'(z)$, and for this we use the recurrence relation given in \cite[Equation 23.9.5]{NIST:DLMF}. The $x$ and $y$ are modular functions for $X_{G}$, but they do not necessarily live in the function field of $X_{G}$ over $\Q$, because the cusp at infinity on $X_{G}$
need not be rational.

\begin{step} \label{Step10} Express the $x$ and $y$ modular functions as ratios of elements in the graded ring $R \otimes \Q(\zeta_{N})$. \end{step}

\noindent To generate the graded ring $R$, we compute the initial ideal of the ideal generated by the equations for $X_{G}$.
Using this, we generate the degree $d$ piece of $R$ (for values of $d$ less than or equal to the upper bound in Step~\ref{Step7}) by enumerating degree $d$ monomials
that are not in the initial ideal, obtaining each as the product of a degree $1$ element and a degree $d-1$ element. Once this is complete, we search for elements $a$, $b$ and $c$ in the degree $d$ piece of $R \otimes \Q(\zeta_{N})$ such that $cx + a = 0$ and $cy + b = 0$ by writing $a,b,c$ as linear combinations of the degree $d$ monomial basis and solving the resulting linear system obtained from the Fourier expansions at the cusps. The map from $X_{G}$ to the elliptic curve is then given by $P \mapsto (-a(P) : -b(P) : c(P))$.
Our goal at this stage is just to evaluate the image of the rational base point $Q$, and so we keep searching for solutions to the linear equations (increasing the degree if necessary) until we find $a$, $b$ and $c$ such that $a(Q), b(Q)$ and $c(Q)$ are not all zero.
We then compute the Fourier expansions of the difference between $(x,y)$ (as computed in Step~\ref{Step9}) and the image of the base point. This gives modular functions $x'$ and $y'$ in the function field of
$X_{G}$ over $\Q$ that satisfy the defining equation for $E'$.

\begin{step} \label{Step11} Express the $x'$ and $y'$ modular functions as ratios of elements in the graded ring $R$. \end{step}

\noindent We use a very similar process as in the last step to find elements $a$, $b$ and $c$ of $R$ so that
$cx' + a = 0$ and $cy' + b = 0$. When we do the linear algebra to find $a$, $b$ and $c$, we frequently find that
the dimension of the null space is greater than $1$. We store up to five different sets of polynomials defining the same map, as this makes
Step~\ref{Step13} more efficient.

\begin{step} \label{Step12} Check the result. \end{step}

\noindent It is necessary to verify that the polynomials obtained in Step~\ref{Step11} do in fact define a map $X_{G} \to E'$ defined over $\Q$. To check this, we could use Magma's built-in map constructor which checks the map by default. However, this is too inefficient when the genus is large so we instead take an alternative approach. We plug the Fourier expansions of the canonical ring elements $a$, $b$, and $c$ into the equation for the elliptic curve $E'$ and determine a lower bound for the order of vanishing of this modular form at each cusp. If the sum of these lower bounds exceeds $\frac{k}{12} [\SL_{2}(\Z) : \Gamma]$, then
by Proposition~\ref{checkingprop}, the elements of the canonical ring obtained by plugging in $a$, $b$ and $c$ into the defining equation for $E'$ are zero.
This proves that the equations we have chosen define a map $X_{G} \to E'$.
We perform this check for each of the different polynomials defining
a map $X_{G} \to E'$ and also check that these polynomials are compatible with one another. The polynomials defining this map are stored in
the {\tt ModEC} record, and the function {\tt FindMapsToEC} returns this record.

\begin{step} \label{Step13} Determine the rational points on $X_{G}$. \end{step}

\noindent This step is handled by the function {\tt RatPtsFromMaps}.
For each point $T \in E'(\Q)$, we compute the preimage scheme $Z = Z_T$ of $T$ inside $X_{G}$ and determine the rational points on this scheme. These schemes can be quite complicated, and for this, rather than using the built-in {\tt RationalPoints} command,
we use a novel method for determining the rational points on these zero-dimensional schemes when possible. 

First, we clear denominators to consider $Z$ as a scheme over $\operatorname{Spec} \Z$. We search for a prime $p$ for which $Z$ has no $\F_{p}$-rational singular points. (Such a prime need not exist.) We use Hensel's lemma to lift each mod $p$ point to a point in $Z(\Z/p^{k} \Z)$ for some $k > 1$. (We take $k = \lceil 12 \log_{p}(10) \rceil + 1$ so that $p^k > 10^{12}$ which allows us to correctly identify rational numbers of height up to about $10^6$.) For each coordinate of the lifted point, we find the simplest rational number that reduces to it mod $p^{k}$ and we test whether this point yields a point in $Z(\Q)$. This process gives us a lower bound on $|Z(\Q)|$.

Next, for primes $\ell > p$ with $\ell \nmid N$, we test whether $Z(\F_{\ell})$ contains singular points.
If it does not, then every point in $Z(\F_{\ell})$ has a unique lift to a point in $Z(\Q_{\ell})$ by Hensel's lemma (see for example \cite[Theorem 3.3]{KConradHensel}). This implies that
$|Z(\Q)| \leq |Z(\Q_{\ell})| = |Z(\F_{\ell})|$. We hope to find a prime $\ell$ for which this upper bound matches the lower bound we found using the prime $p$. In the event that no prime $p$ is found, or for which no prime $\ell < 400$ is found, we fall back to Magma's built-in function {\tt RationalPoints} to determine $Z(\Q)$. After $Z(\Q) = Z_{T}(\Q)$ is determined for all $T \in E'(\Q)$, we have that
\[ X_G(\Q) = \bigcup_{T \in E'(\Q)} Z_T(\Q). \]

\begin{step} \label{Step14} Determine $j(P)$ for $P \in X_{G}(\Q)$. \end{step}

\noindent Once Step~\ref{Step13} is complete, we know $|X_{G}(\Q)|$. We use data from the LMFDB
to determine the number of rational cusps that $X_{G}$ has, and we use
the algorithm of Novak \cite{Novak} to compute the number of rational points on $X_{G}$ with each CM discriminant. It is often the case that the number of rational cusps plus the number of rational CM points is equal to $|X_{G}(\Q)|$; in this case, we have already determined the multiset $\{ j(P) : P \in X_{G}(\Q) \}$ and
no further work is necessary. Otherwise, there exist non-cuspidal, non-CM rational points, and we next check whether the LMFDB stores a model $X$ of our modular curve and a $j$-map $j \colon X \to \P^{1}$. (Note that $X$ and $X_{G}$ are isomorphic over $\Q$, but they need not have the same model.) If such a model and $j$-map are available,
we point search on the LMFDB model, check that we found $|X_{G}(\Q)|$ points on $X$, and compute their images under $j$. If no such model and $j$-map are available, then we could use our intrinsic ${\tt ComputeJ}$
to compute $j \colon X_{G} \to \P^{1}$. During the course of our computation, this was never needed.

\subsection{Discussion of parameters}
Next, we describe a few optional parameters that can be employed. The parameter {\tt Verbose} defaults to false, but when set to \texttt{true} it prints quite a bit of detail about all the steps in the computation. 

In some cases, not enough power series precision is available in Step~\ref{Step3}. Also, it is not known how much power series precision is necessary to perform Step~\ref{Step12} until the degree of the map $X_{G} \to E'$ is known, and this is not computed until Step~\ref{Step7}. The optional parameter {\tt precmult} defaults to $1$, but when set to a number higher than $1$, extra power series precision is used. In the event that not enough precision is available, Step~\ref{Step7} will report a runtime error and suggest rerunning the example with a given value of {\tt precmult}.

Finally, there is an optional parameter {\tt IgnoreBase}. This defaults to false, but when set to \texttt{true}, the code \emph{assumes} that the map to the elliptic curve, with $Q$ chosen to be the point at infinity, is defined over $\Q$. In this situation, Steps~\ref{Step2} and \ref{Step10} are skipped. Examples suggest that this assumption is more likely to hold when the genus of $X_{G}$ is high, and it allows some cases to be handled where $X_{G}(\Q) = \emptyset$ or point searching on $X_{G}$ (using Magma's built-in {\tt PointSearch}) is very slow. If Step~\ref{Step2} does not find a local obstruction
and also does not find a rational point, then the code sets {\tt IgnoreBase} to be true and proceeds (hence skipping Step~\ref{Step10}). In some cases, the assumption that the map $X_{G} \to E$ is defined over $\Q$
is incorrect and the code fails. For more details about these cases,
see Section~\ref{section:unique_cases}.

\section{Illustrative example}
\label{section:ex}

To illustrate the algorithm described in Section~\ref{section:alg}, we trace through each step using the modular curve $X_G$ with label \mc{36.108.6.g.1}. Among all curves to which we applied our method, genus~6 is the highest genus for which a non-cuspidal, non-CM rational point occurs, and there are exactly two such curves (the other is \mc{20.120.6.a.1}). The curve $X_G$ can be realized as the fiber product
\[
X_G = X_{\rm ns}^+(4) \times_{X(1)} \mc{9.27.0.a.1},
\]
where \mc{9.27.0.a.1} is the well-known Elkies curve appearing in \cite{elkies2006ellipticcurves3adicgalois}. The curve $X_G$ appears in \cite[Remark 25]{MR4458130}, although the authors do not determine its rational points. To our knowledge, the rational points on $X_G$ have not previously been computed. 

The curve \mc{36.108.6.g.1} has level $36$, index $108$, genus $6$, and associated group
\[
G = \left\langle \begin{bmatrix}8&17\\25&12\end{bmatrix}, \begin{bmatrix}12&23\\5&15\end{bmatrix}, \begin{bmatrix}14&15\\25&5\end{bmatrix} \right\rangle \leq \GL_2(\Z/36\Z).
\]
From the curve's LMFDB page, we see that its Jacobian decomposes over $\mathbb{Q}$ as
\[
J(X_G) \sim E \times E_1 \times E_2 \times E_3 \times A
\]
where $E$ is an elliptic curve of rank $0$, $E_1$, $E_2$, and $E_3$ are distinct elliptic curves of rank $1$, and $A$ is an abelian surface. The elliptic curve $E$ is isogenous to \href{https://www.lmfdb.org/EllipticCurve/Q/432/f/1}{\texttt{432.f1}}, which is given by the equation
\[
E \colon y^2=x^3-27x-918.
\]
In general, our code constructs a map $X_G \to E'$, where $E'$ is some curve in the isogeny class of $E$. In the present example, the isogeny class of $E$ has size $1$, so necessarily $E' = E$.

In the previous paragraph, we used the Jacobian decomposition of $X_G$. Alternatively, one could use our intrinsic \texttt{EllipticCurveQuoCandidates}, which produces a finite superset of elliptic curves from Cremona's database \cite{MR2282912} that could appear in the Jacobian decomposition. By taking the precision sufficiently high, this superset is exactly the set of elliptic curves from Cremona's database appearing in the Jacobian decomposition. Note that an elliptic curve appearing in the Jacobian decomposition of $X_G$ must have conductor dividing the square of the level of $X_G$.

Our \texttt{Magma} implementation first transposes  each generator of $G$ to obtain generators for the subgroup $G^\intercal$,
\[
G^\intercal = \left\langle \begin{bmatrix}8&25\\17&12\end{bmatrix}, \begin{bmatrix}12&5\\23&15\end{bmatrix}, \begin{bmatrix}14&25\\15&5\end{bmatrix} \right\rangle \leq \GL_2(\Z/36\Z).
\]
Calling Zywina's code on generators of $G^\intercal$ determines that $X_{G}$ is not geometrically hyperelliptic and produces the canonical model of $X_G \subset \mathbb{P}^5$. 
This model is obtained from the $6$ dimensional $\Q$-vector space of weight $2$ cusp forms for $G$. The code then selects the working power series precision. In this example, the  precision is set to 136.

Since the genus of $X_G$ is less than $30$,  the code performs
a point search on the canonical model using the \texttt{PointSearch} command. The first rational point found is 
\[ Q = (1 : 2 : -2 : 7 : 3 : 4) \in X_G(\Q).\]
This point is used as the base point for the Abel--Jacobi map and later when translating the map to ensure it is defined over $\Q$. Because a rational point was found, the local solvability checks are skipped.

Next, we seek a weight $2$ cusp form whose Hecke eigenvalues match those of $E$ away from $N = 36$.  Concretely, we iterate over primes $p \nmid 36$ and intersect the kernels of the operators
\[ T(p) \begin{bmatrix} 1 & 0 \\ 0 & p^{-1} \end{bmatrix}^{*} - a_{p}(E). \]
Taking $p = 5$, the kernel is already $1$-dimensional.  Choosing a
generator and scaling it to have coprime algebraic integer coefficients, we obtain the desired eigenform $f$.

Next, we compute generators of $G \cap \SL_{2}(\Z/N\Z)$ and lift them to $\SL_2(\Z)$,  obtaining
\[
\left\{
\begin{bmatrix}
17 & 144 \\
36 & 305
\end{bmatrix},
\begin{bmatrix}
23 & 394 \\
15 & 257
\end{bmatrix},
\begin{bmatrix}
12 & -53 \\
17 & -75
\end{bmatrix},
\begin{bmatrix}
31 & 143 \\
13 & 60
\end{bmatrix},
\begin{bmatrix}
19 & -180 \\
36 & -341
\end{bmatrix} \right\}.
\]
We then adjoin additional random matrices in the kernel of the reduction map 
\[ \SL_{2}(\Z) \to \SL_{2}(\Z/N\Z),\] 
until we have $20$ matrices in total. Each matrix is expressed as a word in $S = \begin{bsmallmatrix} 0 & -1 \\ 1 & 0 \end{bsmallmatrix}$ and $T = \begin{bsmallmatrix} 1 & 1 \\ 0 & 1 \end{bsmallmatrix}$, up to sign. For example,
\[
\begin{bmatrix}
17 & 144 \\
36 & 305
\end{bmatrix} = -
S^{1} T^{-3} S^{1} T^{-2} S^{1} T^{-2} S^{1} T^{-2} S^{1} T^{-2} S^{1} T^{-2} S^{1} T^{-2} S^{1} T^{-2} S^{1} T^{-3} S^{1} T^{8}.
\]
By default, the random matrices are chosen so that these words have length $\leq 50$. For each of the $20$ matrices $g$, we compute the complex period
\[
\phi(g) = 2 \pi i \int_{i}^{g(i)} f(s) \, ds
\]
by breaking the path into segments corresponding to words in $S$ and $T$ and repeatedly applying \eqref{E:transformation}. After rounding to two decimal places for brevity, the computed set of periods is
{\scriptsize \[
\left\{
\begin{array}{c}
326.76, -326.76, 163.38, 0, 81.69 + 150.04i, -81.69 + 150.04i, -81.69 - 150.04i, \\
81.69 - 150.04i, 408.45 + 150.04i, -408.45 + 150.04i, -326.76 + 300.08i, -81.69 + 450.12i, \\ -571.83 + 150.04i, 163.38 + 300.08i, -163.38 - 300.08i, 163.38 - 300.08i
\end{array}
\right\}.
\]}

An elliptic curve over $\Q$ has a period lattice of the form
\[ \{ a \omega_{1} + b \omega_{2} : a, b \in \Z \} \]
where $\omega_{1}$ is real and either $\omega_{2}$ is purely imaginary or satisfies ${\rm Re}(\omega_{2}) = \frac{1}{2} \omega_{1}$. From the numerical output above, we recognize
the second shape and take
\[
\omega_1 \approx 163.379 \quad \text{and} \quad \omega_2 \approx 81.69 + 150.039i
\]
so that the computed periods are consistent with the $\mathbb{Z}$-span of $\omega_1$ and $\omega_2$. At this stage, it is in principle possible that our finite sampling reveals only a sublattice of the true period lattice. 

In the next step, we compare the above lattice to the period lattices of each curve in the isogeny class of $E$. In this example, the isogeny class of $E$ has size $1$. A $\mathbb{Z}$-basis for the period lattice of $E$ is $\omega_1' \approx 0.756$ and $\omega_2' \approx 0.378 + 0.695i$.
Comparing the ratios gives
\(
\frac{\omega_1'}{\omega_1} \approx \frac{\omega_2'}{\omega_2} \approx \frac{1}{216}.
\)
After scaling, the computed period lattice matches closely with that of $E$. Thus, the Manin constant is heuristically $\frac{1}{216}$, associated with the optimal curve $E' = E$. 

Using \cite[Theorem 3]{Cremona}, we compute that the analytic modular degree of the map $X_G \to E$ is $6$. In this example,  this computation involves evaluating $108$ periods and is the slowest part of the overall computation, taking about 19 seconds.

The cusps of $X_G$ are $\infty$, $\frac{5}{36}$, and $\frac{7}{36}$. For each cusp $P$, we numerically compute the integral
\( \int_{\infty}^{P} \tfrac{1}{216} f(s) \, ds, \)
which we find lies in $\Lambda(E)$ in each case. The code expresses each as a linear combination of the basis periods and recovers the exact torsion point in $E(\Q(\zeta_{36}))$ by factoring primitive division polynomials over $\Q(\zeta_{36})$.

We next compute the Fourier expansions of the modular functions $x$ and $y$ defining the map $X_G \to E$ at all cusps using the process described in Step~\ref{Step9}. Because the modular degree is $6$, we know that we can realize $x$ and $y$ as ratios of degree $3$ elements of the graded ring $R \otimes \Q(\zeta_{36})$, which allows us to reduce the precision from $136$ to $64$.

Next, we realize $x$ and $y$ in the function field over $\Q(\zeta_{36})$. To do so, we compute a Gr\"obner basis for the ideal defining the model of $X_G$ and use the resulting initial ideal to enumerate a basis of monomials in each degree. We set up the homogeneous linear system encoding $cx+a=0$ and $cy+b=0$ at all cusps and take its null space. This yields degree 3 triples $(a,b,c)$ and hence a map $P\mapsto(-a(P):-b(P):c(P))$; we choose one with $(a(Q),b(Q),c(Q))\neq(0,0,0)$ where $Q$ is our chosen base point.

We now have a map $X_G \to E$ defined over $\Q(\zeta_{36})$. We evaluate the image of the rational base point $Q$ in $E(\Q(\zeta_{36}))$ and translate so that $Q$ maps to the origin. Subsequently, we repeat the linear algebra over $\Q$ rather than $\Q(\zeta_{36})$, producing one or more sets of homogeneous polynomials with rational coefficients defining the same map. (In this example, the image of $Q$ is the point at infinity on $E$ so the translation is trivial; this shows in particular that the example would run  successfully with the parameter \texttt{IgnoreBase} set to \texttt{true}.)

We now have a map $\psi \colon X_G \to E$ defined over $\Q$ given by six homogeneous cubic polynomials. 
We must check that this is indeed a nonconstant morphism whose image lies in the intended codomain.
We use our function \texttt{MapOverQ} to quickly certify correctness using Proposition~\ref{checkingprop}.

The Mordell--Weil group of $E$ is
\( E(\Q) = \{ (0:1:0) \}. \)
Pulling back $(0:1:0)$ along $\psi$, we find that
\[
X_G(\Q) = \{ (0 : 1 : -1 : -2 : 0 : 1), (1 : 2 : -2 : 7 : 3 : 4), (2 : 1 : -1 : 2 : 6 : 5) \}.
\]
This can be accomplished using our intrinsic \texttt{RatPtsFromMaps} or, because the genus is relatively low, by directly taking the preimage of $(0:1:0)$ in Magma, which gives a zero dimensional scheme whose rational points can be determined by calling \texttt{RationalPoints}.

In the previous step, we determined that $X_G$ has exactly three rational points. The LMFDB page for this modular curve includes the canonical model and the $j$-map from this model to $\P^1$. To avoid recomputing the $j$-map, which can be done using our intrinsic \texttt{ComputeJ} but is time intensive, we instead perform a point search on the LMFDB's model to find the three rational points and evaluate their images under the provided $j$-map. Doing so, we find that two of the points have the CM $j$-invariant $j = 0$ and one has the non-CM $j$-invariant $j = -2^2 \cdot 3^7 \cdot 5^3 \cdot 439^3$. The total time for handling this case was 31.58 seconds. If we do not use the LMFDB's model and compute the $j$-map using \texttt{ComputeJ}, an additional 1289.87 seconds are required.

\section{Computational results}
\label{section:results}

There are 1034 minimal modular curves in the sense of Section~\ref{section:intro}. In this section, we describe our computational run on these curves, note the cases that presented difficulties, and discuss the performance we observed. Since the runtime of our code generally increases with the genus, we structured our computations according to the genus of the modular curves.

We begin with the minimal modular curves of genus $1 \leq g \leq 10$. There are $666$ such curves. Based on data from the LMFDB (which uses \cite[Section 5]{RSZB}), we know that $91$ of them have a $p$-adic obstruction for some prime not dividing the level. In addition, using our \texttt{SmartpAdic} function, we find that at least $117$ have a $p$-adic obstruction for some prime dividing the level. For each of the remaining $458$ curves, we ran our main intrinsic \texttt{FindMapsToEC} and, when successful, the intrinsic \texttt{RatPtsFromMaps}.

For the $230$ minimal modular curves of genus $11 \leq g \leq 60$, we are more selective and do not run our code when we observe any of the following situations:
\begin{enumerate}
    \item \label{Sit1} There is a $p$-adic obstruction for some prime $p$ not dividing the level.
    \item \label{Sit2} The genus is $\leq 29$ and there exists a prime $p$ dividing the level for which the curve has no points modulo $p^k$ (where $k$ is as in Section~\ref{section:alg}).
    \item \label{Sit3} The curve covers a prime power level modular curve with finitely many rational points, all of which are already known from prior work \cite{Cursed,MoreQuadChab,RSZB}.
    \item \label{Sit4} There exist distinct primes $p \geq 5$ and $q$ such that the mod $p$ reduction of $G$ is contained in the normalizer of the non-split Cartan subgroup modulo $p$ and the mod $q$ reduction of $G$ is contained in the Borel subgroup or in the normalizer of the split Cartan subgroup modulo $q$. Here there are no non-CM rational points by \cite{Lemos1,Lemos2}; see also \cite[Theorem 6.2]{DanielsRouse} and \cite[Lemma 17]{2512.00652}.
\end{enumerate}
In any of these situations, the rational points are known and we do not run our code. After applying these criteria, 38 minimal curves in the genus range $11 \leq g \leq 60$ remained, and for each of these we ran \texttt{FindMapsToEC} and, when successful, \texttt{RatPtsFromMaps}.

There are $138$ minimal modular curves of genus $\geq61$. After applying the four criteria above, the rational points are already known for $127$ of these curves. Further, we handle one of the remaining eleven curves, \mc{63.1512.115.t.1}, in Section \ref{section:unique_cases}. The other ten curves are precisely those that we were unable to handle in Theorem~\ref{theorem:main}. As an illustration of the computational difficulty in this range, we attempted the case of \mc{62.930.63.a.1} where the degree of the map $X_{G} \to E'$ is $240$. We stopped the Magma process after more than one month had been spent in Step~11 trying to compute the null space of a very large matrix.

Altogether, we ran our main code on 496 minimal curves (excluding the  genus $\geq 61$ curves that appear out of reach). With the default parameters, the computation succeeded in all but 27 cases:
\begin{itemize}
    \item In 16 cases, the code failed in Step~\ref{Step3} to identify the space of Hecke eigenforms because insufficient power series precision was available.
    \item In 7 cases, the code failed in Step~\ref{Step7} after it was determined that there would not be enough precision to check the map.
    \item In 3 cases (\mc{20.60.2.j.1}, \mc{52.182.10.a.1}, and \mc{56.126.6.c.1}), no rational base point was found and the code proceeded with the parameter {\tt IgnoreBase} set to \texttt{true}. When this happens, the code assumes that the map $X_{G} \to E'$ obtained using the cusp at infinity as the base point is defined over $\Q$; in each of these three cases, that assumption is false.
    \item Finally, there is one case (\mc{42.168.11.c.1}) in which the code failed to match the period lattice of one of the six modular forms $f$ with the period lattice of an elliptic curve in the isogeny class of $E$.
\end{itemize}
We handle \mc{20.60.2.j.1}, \mc{52.182.10.a.1}, and \mc{56.126.6.c.1} in Section~\ref{section:unique_cases}. For the remaining 24 failures, we reran the computations with increased precision. All reruns were successful except for \mc{66.660.49.e.1} and \mc{66.660.50.a.1}. For \mc{66.660.49.e.1}, a segmentation fault occurred when Step~\ref{Step13} was unable to use Hensel lifting to determine $X(\Q)$ and fell back to Magma's built-in {\tt RationalPoints}. Rerunning this case led to success. For \mc{66.660.50.a.1}, there was still insufficient precision to compute the Hecke eigenform, so we reran with {\tt precmult} set to $4.5$. In this instance, a segmentation fault again occurred in Step~\ref{Step13}, and rerunning the script several times led to the same outcome. This suggests that the map $X_{G} \to E'$ is ramified at the point at infinity, which is the only rational point of $E'$. We likewise handle \mc{66.660.50.a.1} in Section~\ref{section:unique_cases}.

We next briefly discuss the degrees of the maps $X_{G} \to E'$ that were found. Regarding small degrees, the highest genus modular curves we found admitting a degree $2$ map to an elliptic curve have genus $7$. There are two such curves: $X_{0}(62)$ and $X_{0}(69)$, which were shown to be bielliptic in \cite{Bars}. At the other extreme, the highest degree maps to elliptic curves that we computed have degree $126$, occurring for \mc{42.126.10.b.1} and \mc{42.168.12.b.1}. These were also the two most time-consuming examples we ran, requiring roughly 12.5 and 9.5 hours, respectively. In the first case, finding the map $X_{G} \to E'$ requires realizing $x'$ and $y'$ as ratios of degree $22$ elements of the canonical ring $R$. 

Our code computes the modular degree for each candidate rank $0$ isogeny class. This is costly, because computing the modular degree is often the slowest part of the computation. For some of the modular curves we encountered, there are as many as nine different isogeny classes of rank $0$ elliptic curves to consider, and the modular degree computation is carried out for each such class. However, this approach allows us to select a map of minimal degree, which can lead to time savings. For example, for the genus $54$ curve \mc{63.756.54.a.1}, there are two isogeny classes of elliptic curves of rank $0$ in the Jacobian decomposition. For one of these ($y^{2} + y = x^{3} - 27x - 7$) the degree of the map is $72$, while for the other ($y^{2} + y = x^{3} - 6x + 3$), the degree of the map is $216$.

The total computation time across all runs was about 79 CPU hours. The most time-consuming steps in the algorithm are computing the modular degrees, expressing $x$ and $y$
as ratios of elements of $R \otimes \Q(\zeta_{N})$, expressing $x'$ and $y'$ as ratios of elements of $R$, and pulling back the rational points on $E'$. These steps account for approximately 45\%, 8.7\%, 22\%, and 9.4\% of the total computation time, respectively. For each genus $g \leq 9$, the average runtime per curve of genus $g$ was approximately one minute or less. For higher genera, the runtimes we observed were more variable and generally higher.

After this paper had been submitted, the authors noticed two flaws in the code. First, there was an inefficiency in how generators for the degree $d$ parts of the ring $R$ were chosen that was unnecessarily slow and memory intensive. Second, although multiple defining polynomials for the map $X_{G} \to E'$ were being computed, only one was being used in {\tt RatPtsFromMaps}. We corrected these and reran all the computations; a second set of log files is available at \cite{ModCrvToEC}. With these changes, determining rational points on \mc{66.660.50.a.1} was now successful. 

Finally, we note that of the $1024 (=1034-10)$ minimal modular curves that we handle in this article, only $33$ have a non-CM rational point. Non-CM rational points on $X_{G}$ only occur
when the genus is $\leq 6$. Of these $33$ minimal modular curves with a non-CM rational point, $18$ have genus $1$, $10$ have genus $2$, $2$ have genus $3$, $1$ has genus $4$, and $2$ have genus $6$. All of the non-CM $j$-invariants 
are listed in Table~\ref{jinvtable}.

\section{Unique cases}
\label{section:unique_cases}

In this section, we complete the determination of rational points on the five remaining modular curves that were not resolved in our main computational run. Three of the curves (\mc{20.60.2.j.1}, \mc{52.182.10.a.1}, and \mc{56.126.6.c.1}) have a possible local-global failure and two (\mc{66.660.50.a.1} and \mc{63.1512.115.t.1}) have high genus.

We were unsuccessful in determining the rational points on \mc{20.60.2.j.1},\\ \mc{52.182.10.a.1}, and \mc{56.126.6.c.1} because we found no rational points through point searching, yet observed no local obstruction, and our code was unsuccessful with \texttt{IgnoreBase} set to \texttt{true}. In addition to these three curves, we encountered the curves \mc{26.182.10.a.1} and \mc{60.120.7.da.1}, which may also have a local-global obstruction, but for which our code ran successfully with \texttt{IgnoreBase} set to \texttt{true}.

The curve \mc{20.60.2.j.1} is genus $2$ rank $0$ and has a model
\[
y^2 = - x^6 - 2x^5 - 5x^4 + 5x^3 + 5x^2 + 18x + 11.
\]
Running \texttt{Chabauty0}, we find that this curve has no rational points.

The curve \mc{52.182.10.a.1} covers the curve $X_{S_4}(13)$, which has only three non-CM rational points \cite[Section 5.1]{MoreQuadChab}. The associated $j$-invariants are
\scriptsize \[
j_1 = \frac{2^4 \cdot 5 \cdot 13^4 \cdot 17^3}{3^{13}}, \;
j_2 = -\frac{2^{12} \cdot 5^3 \cdot 11 \cdot 13^4}{3^{13}}, \; \text{and} \;
j_3 = \frac{2^{18} \cdot 3^3 \cdot 13^4 \cdot 127^3 \cdot 139^3 \cdot 157^3 \cdot 283^3 \cdot 929}{5^{13} \cdot 61^{13}}.
\] \normalsize
Using \texttt{Magma}, we construct elliptic curves $E_{j_1}$, $E_{j_2}$, and $E_{j_3}$ with these $j$-invariants. Using Zywina's \texttt{FindOpenImage} \cite{ZywinaOpen}, we compute their adelic images and observe that, for each of the three curves, the mod $52$ reduction is not contained in the group associated with \mc{52.182.10.a.1}. Consequently, none of the three non-CM rational points on $X_{S_4}(13)$ lift to \mc{52.182.10.a.1}, and hence this curve has no non-CM rational points. Since \mc{52.182.10.a.1} has no rational cusps or CM points, it has no rational points.

The curve \mc{56.126.6.c.1} covers \mc{14.63.2.a.1}, which is genus $2$ rank $0$ and has a model
\[
y^2 = -8x^5 + 17x^4 + 2x^3 - 11x^2 + 4x + 4.
\]
Running \texttt{Chabauty0}, we find that this curve has exactly four rational points, all of which are CM. Since \mc{56.126.6.c.1} has no CM points, it follows that it has no rational points.

The modular curve \mc{66.660.50.a.1} has genus $50$ and presented difficulties for our initial computation, so we treat it separately here. It covers the curve \mc{33.110.4.a.1}, whose canonical model is given on its LMFDB page. Let $X_{{\rm ns}}^{(-3)}(11)$ denote the quadratic twist by $-3$ of the double cover $X_{{\rm ns}}(11) \to X_{{\rm ns}}^+(11)$ given by 
\[
\begin{cases}
    y^2 = 4x^3 - 4x^2 - 28x + 41 \\
    t^2 = 3(4x^3 + 7x^2 - 6x + 19).
\end{cases}
\]
Since \(\mc{33.110.4.a.1}\) and \(X_{{\rm ns}}^{(-3)}(11)\) have the same level and genus, and both admit a degree \(2\) map to \(X_{\rm ns}^{+}(11)\), it is plausible that these curves are isomorphic. Indeed, using Magma, we verify that these two curves are isomorphic over $\Q$. In \cite[Section 7.1]{MR4036449}, the authors consider \(X_{{\rm ns}}^{(-3)}(11)\) and determine its rational points, which are $(-2, \pm 7, \pm 9)$ and $(4, \pm 11, \pm 33)$. Thus \mc{33.110.4.a.1} has exactly $8$ rational points. Exactly two of these are non-CM, both with $j$-invariant
\[
j = \frac{2^{18} \cdot 3^3 \cdot 5^3 \cdot 7^1 \cdot 11^3 \cdot 23^3 \cdot 29^3 \cdot 103^3}{67^{11}}.
\]
We construct an elliptic curve $E_{j}$ with this $j$-invariant. Again using Zywina's \texttt{FindOpenImage}, we compute the adelic image of $E_j$ and observe that its mod $66$ reduction is not contained in the group associated with \mc{66.660.50.a.1}. Consequently, neither of the two non-CM rational points on \mc{33.110.4.a.1} lift to \mc{66.660.50.a.1}, and hence this curve has no non-CM rational points.

During our computation, we did not attempt to determine the rational points on \mc{63.1512.115.t.1} due to its high genus, and we instead handle it here individually. This curve covers \mc{21.84.3.a.1}, which we denote by $X$. This curve does not admit a map to a rank $0$ elliptic curve defined over $\Q$, so this curve is not handled as part of our main run. We note that $X$ may be of independent interest because
it is the fiber product $X_{{\rm ns}}^{+}(3) \times_{X(1)} X_{{\rm sp}}^{+}(7)$. The remainder of this section is devoted to determining the rational points of $X$. 

The Jacobian of $X/\Q$ is the product of a rank $1$ elliptic curve and a rank $0$ abelian surface $A/\Q$ corresponding to the newform orbit with LMFDB label \href{https://www.lmfdb.org/ModularForm/GL2/Q/holomorphic/63/2/a/b/}{\tt 63.2.a.b}. On the page for this newform orbit, the LMFDB notes that $A$ is isogenous to the Jacobian of the genus $2$ curve
with label \mc{3969.d.250047.1},
\[
y^{2} + (x^{2} + x + 1)y = -3x^{5} + 5x^{4} - 4x^{3} + x.
\]
The page for this curve notes that its Jacobian factors as the square of an elliptic curve of rank $0$ (isogeny class \href{https://www.lmfdb.org/EllipticCurve/2.0.3.1/441.2/a/}{\tt 2.0.3.1-441.2-a}) over $\Q(\sqrt{-3})$. Hence there is a map from $X$ to an elliptic curve over $\Q(\sqrt{-3})$ that has rank $0$. We determine $X(\Q)$ by writing down this map, following as closely as possible the algorithm in Section~\ref{section:alg}. 

The main
challenge is to identify an appropriate weight $2$ cusp form that corresponds to a $1$-form on $X/\Q(\sqrt{-3})$ whose periods
match those of one of the eight elliptic curves in the isogeny class \href{https://www.lmfdb.org/EllipticCurve/2.0.3.1/441.2/a/}{\tt 2.0.3.1-441.2-a}. We first compute the action
of $T(2)$ on the 3-dimensional space of weight $2$ cusp forms. This Hecke operator has a $1$-dimensional kernel (corresponding
to the map from $X$ to a rank $1$ elliptic curve over $\Q$), and the image of the Hecke operator is a $2$-dimensional Hecke-invariant subspace. We choose a basis $\{ f_{1}, f_{2} \}$ for this space and compute periods of $f_{1}$ and $f_{2}$ for $30$ matrices. Numerically, we observe that
there are two real numbers $x_{1} \approx 11.94808$ and $x_{2} \approx 11.30988$ such that
every period of $f_{1}$ or $f_{2}$ is an integer linear combination of
$x_{1}$, $x_{1} \sqrt{-3}/3$, $x_{2}$, and $x_{2} \sqrt{-3}/3$. 

Next, we take one elliptic curve
in the isogeny class {\tt 2.0.3.1-441.2-a} and express its periods $\Omega_{1}$ and $\Omega_{2}$ in terms of $x_{1}$ and $x_{2}$. We then seek a linear combination of $f_{1}$ and $f_{2}$ with
coefficients in $\Q(\sqrt{-3})$ whose periods all lie in the $\Q$-span of $\Omega_{1}$ and $\Omega_{2}$. The form $f_{1} + f_{2} \sqrt{-3}$ works, and the periods of
$\frac{1}{21} (f_{1} + f_{2} \sqrt{-3})$ exactly match those of the elliptic curve \href{https://www.lmfdb.org/EllipticCurve/2.0.3.1/441.2/a/3}{\tt 441.2-a3}
\[
E : y^{2} + xy + ay = x^{3} - x^{2} + (4a-3)x - 4a + 1
\]
where $a^{2} - a + 1 = 0$. Using this information, we can proceed with the rest of the algorithm.

It produces a degree $6$ map $\phi \colon X \to E$ defined over $\Q(\sqrt{-3})$. The Mordell--Weil group of $E$ is $E(\Q(\sqrt{-3})) \simeq \Z/6\Z$ and pulling back these points we find that
$X(\Q(\sqrt{-3}))$ consists of five points, of which three are rational. The three rational points are a rational cusp, a CM point with discriminant $-3$,
and a CM point with discriminant $-19$. The two non-rational points in $X(\Q(\sqrt{-3}))$ are also CM points with discriminant $-19$. Details of these computations are in the file {\tt level21genus3map.m}.  

\section{The remaining cases}
\label{section:remaining_cases}

Combining our computations of Section~\ref{section:results} with the work of Section~\ref{section:unique_cases}, we have determined the rational points on all $1034$ minimal modular curves, with the exception of the ten modular curves listed in Theorem~\ref{theorem:main}. To complete the proof of Theorem~\ref{theorem:main}, we also need to consider modular curves that minimally cover one of these ten curves and have level $\leq 70$. All of these curves fall into either situation \eqref{Sit1}, \eqref{Sit3} or \eqref{Sit4} from Section~\ref{section:results} and so they have no non-CM rational points.

For example, let us consider the modular curve \mc{55.3300.239.a.1}. There are five curves of level $\leq 70$ that cover this curve, none of which have non-CM rational points, as we now explain. The curves \mc{55.6600.489.a.1} and \mc{55.6600.489.b.1} have no real points and \mc{55.6600.489.d.1} has no points over $\Q_2$. The curve \mc{55.6600.489.c.1} covers \mc{11.110.4.b.1}, which admits a map to an elliptic curve over $\Q$ of rank $0$ and is handled as part of the computations of this article. It has no non-CM rational points, which was first determined by Zywina \cite[Lemma 4.5]{Zywina2015}. Lastly, \mc{55.9900.730.d.1} covers $X_{S_4}(11)$, which is an elliptic curve over $\Q$ of rank $0$ and has no non-CM rational points; this was first proved in \cite[Proposition 4.4.8.1]{Ligozat}.

\subsection*{Competing interest statement} 

The authors have no competing interests to declare that are relevant
to the content of this article.

\subsection*{Data availability and reproducibility}

The code accompanying this paper is available at the GitHub repository \cite{ModCrvToEC}. The computations for this article were done using Magma \cite{Magma} V2.29-10 on a desktop with
an Intel Core i9-11900K CPU with 128 GB of RAM, running Ubuntu 24.04. The total run time was about 79 CPU hours and up to 19 GB of RAM was required at some points. The log files from the computations are available at the \href{https://github.com/rouseja/ModCrvToEC/tree/main/logs}{logs} subfolder of \cite{ModCrvToEC}. More information about how to reproduce these results can be found in the file \href{https://github.com/rouseja/ModCrvToEC/blob/main/README.md}{README.md}.

\end{document}